\input amstex
\documentstyle {amsppt}
\magnification=1200
\vsize=9.5truein
\hsize=6.5truein
\nopagenumbers
\nologo

\def\norm#1{\left\Vert #1 \right\Vert}
\def\abs#1{\left\vert #1 \right\vert}

\topmatter

\title
On simplicity of reduced C$^*$-algebras of groups
\endtitle

\author
Pierre de la Harpe
\endauthor

\abstract
A countable group is C$^*$-simple if its reeduced C$^*$-algebra is a simple algebra.
Since Powers recognised in 1975 that non-abelian free groups are C$^*$-simple,
large classes of groups which appear naturally in geometry have been identified,
including non-elementary Gromov hyperbolic groups and lattices in semisimple groups.

In this exposition, C$^*$-simplicity for countable groups 
is shown to be an extreme case of non-amenability.
The basic examples are described and several open problems are formulated. 

\medskip

\noindent
1.~Introduction, and a first definition of \lq\lq C$^*$-simple groups\rq\rq .

\noindent
2.~Formulation in terms of C$^*$-algebras.

\noindent
3.~Powers result on non-abelian free groups.

\noindent
4.~Other examples of C$^*$-simple groups.

\noindent
Appendix.

I.~On the definition of weak containment (Definition~1).

II.~On weak equivalence not preserving irreducibility.

III.~On diagonal and arbitrary matrix coefficients.

IV.~On a definition of Zimmer.

V.~Other notions of equivalence.

VI.~The dual of a group.

VII.~Groups of Type I.

VIII.~The reduced dual.

IX.~The Dixmier property.

X.~On the icc condition.
\endabstract


 \endkeywords

\thanks
Part of this report is on work joint with M. Bekka, 
for which support from the Swiss National Science Foundation 
is gratefully acknowledged.
\endthanks

\address
Pierre de la Harpe, Section de Math\'ematiques, Universit\'e de Gen\`eve,
C.P. 64, CH-1211 Gen\`eve 4, Suisse. 
E-mail : Pierre.delaHarpe\@math.unige.ch
\endaddress

\endtopmatter

\document

\head
{\bf
1. Introduction, and a first definition of \lq\lq C$^*$-simple groups\rq\rq 
}
\endhead

   The subject of this exposition
is that of unitary group representations in Hilbert spaces
and related C$^*$-algebras. 
It is an updated version of notes for a lecture presented at
Bar-Ilan and Jerusalem in June 2000.

\medskip

   Let $G$ be a topological group.
A {\it unitary representation} of $G$ in a Hilbert space $\Cal H$
is a homomorphism 
$\rho : G \longrightarrow \Cal U (\Cal H)$ 
in the unitary group of $\Cal H$;
it is part of the definition that the companion mapping 
$G \times \Cal H \longrightarrow \Cal H$ is continuous. 
Recall that $\Cal U (\Cal H)$ is the group of all linear operators
$u : \Cal H \longrightarrow \Cal H$ which are unitary,
namely isometric ($\norm{u(\xi)} = \norm{\xi}$ for all $\xi \in \Cal H$)
and onto. 

There are two unitary representations which play a leading part.
The first is the {\it unit representation} $1_G$, 
for which $\Cal H$ is the complex field, 
defined by $1_G(g) = \operatorname{id}$ for all $g \in G$.
The second is the {\it left-regular representation} $\lambda_G$,
which is defined only when the group $G$ is {\it locally compact}
(which is by far the most important case);
it acts on the Hilbert space $L^2(G)$ 
of complex-valued measurable functions $\xi$ on $G$
which are square-integrable with respect to a left-invariant Haar measure
(functions modulo equality almost everywhere),
and it is defined by $\left(\lambda_G(g)\xi\right)(g') = \xi(g^{-1}g')$
for all $g,g' \in G$.

The notions of unitary equivalence and of subrepresentations
are well suited for comparing unitary representations of compact groups.
Other groups require weaker and less rigid notions,
defined in terms of coefficients. 
For a representation $\rho : G \longrightarrow \Cal U (\Cal H)$ as above,
we associate to each vector $\xi \in \Cal H$ the corresponding
{\it diagonal matrix coefficient}
$$
\varphi_{\xi,\xi} \, : \, \left\{
\aligned G \ &\longrightarrow \hskip1truecm \bold C \\
         g \ &\longmapsto \ \langle \xi \mid \rho(g) \xi \rangle
\endaligned \right.
$$
which is  a concrete embodiment, namely a continuous complex-valued function, 
of  more abstract data, namely of the representation.  

\medskip

\proclaim{1.~Definition}
Given two unitary representations $\rho,\sigma$ of the same group $G$, 
say that $\rho$ is {\rm weakly contained} in $\sigma$, 
and write $\rho \prec \sigma$, 
if any diagonal matrix coefficient of $\rho$ is 
a limit of sums of diagonal matrix coefficients of $\sigma$, 
uniformly on compact subsets of $G$. 
Two unitary representations are {\rm weakly equivalent,} 
written $\rho \sim \sigma$,  
if each one is weakly contained in the other. 
\endproclaim

\medskip

   A straightforward example is that of a
subrepresentation $\rho$ of a representation~$\sigma$~: 
it is then obvious that $\rho \prec \sigma$. 
But the converse does not hold;
indeed, it is easy to check 
\footnote{
Diagonal matrix coefficients of $1_{\bold R}$ are constant functions.
For $n \ge 1$, define a unit vector $\xi_n \in L^2(\bold R)$ by
$\xi_n(t) = 1/\sqrt{2n}$ if $\vert t \vert \le n$ and
$\xi_n(t) = 0$ if $\vert t \vert > n$;
for any constant $c > 0$, we have
$\lim_{n \to \infty} \left( \sup_{\vert s \vert \le c}
\left\vert 1 - \langle \xi_n \mid \lambda_{\bold R}(s)\xi_n \rangle \right\vert
\right) = 0$,
and it follows that $1_{\bold R} \prec \lambda_{\bold R}$.
On the other hand, $1_{\bold R}$ is not a subrepresentation of $\lambda_{\bold R}$
since the only constant function in $L^2(\bold R)$ is zero.
For any locally compact group~$G$ which is not compact,
the same argument shows that 
$1_{G} \prec \lambda_{G}$
and that $1_{G}$ is not a subrepresentation of~$\lambda_{G}$.
}
that the one-dimensional representation $1_{\bold R}$ of $\bold R$ 
is weakly contained in the regular representation~$\lambda_{\bold R}$, 
and that $1_{\bold R}$ is not a subrepresentation of $\lambda_{\bold R}$. 
In the particular case of two finite-dimensional representations 
$\rho$ and $\sigma$, 
we have $\rho \prec \sigma$ if and only if $\rho$ is a subrepresentation 
of $N \sigma$ (in the naive sense), where $N \sigma$
stands for a direct sum of $N$ copies of $\sigma$, and where
the integer $N$ is large enough; if $\rho$ is moreover irreducible, 
then $\rho \prec \sigma$ if and only  if  $\rho$ is a subrepresentation 
of $\sigma$. \par

See Parts I to III of the appendix for more comments,
and Part IV for a slightly different notion.
See Part V for {\it other} notions of equivalence.

\medskip

   We define now \lq\lq C$^*$-simple groups\rq\rq \
in terms of the notions introduced so far. 
The reformulation of the next section (Corollary~8)
will justify the terminology
and allow us to sketch a general setting.
Also, we recall one of the many standard equivalent definitions
of \lq\lq amenable groups\rq\rq .

\medskip

\proclaim{2.~Definition} A topological group $G$ is 
{\rm amenable} if $1_{G} \prec \lambda_{G}$, 
and {\rm C$^*$-simple} if,
for any unitary representation $\rho$ of $G$, 
the conditions $\rho \prec \lambda_{G}$ and $\rho \sim \lambda_{G}$
are equivalent.
\endproclaim

\medskip

   On the one hand, standard examples of amenable groups include 
soluble groups, locally finite groups, 
and finitely-generated groups of subexponential growth; 
for other examples, see  \cite{CeGH--99}, \cite{GrZu--02}, and \cite{BaVi}.
If $G$ is amenable, 
it is a fact that {\it any} unitary representation of $G$ 
is weakly contained in the left-regular representation~$\lambda_G$
(Proposition 18.3.6 in \cite{DC$^*$--69}). 
On the other hand, C$^*$-simple groups include 
\roster
\item"$\circ$" 
   non-abelian free groups \cite{Pow--75}, 
\item"$\circ$"
   non-trivial free products \cite{PaSa--79}, 
\item"$\circ$"
   non-soluble subgroups of $PSL_2(\bold R)$ \cite{Har--85}, 
\item"$\circ$"
   torsion-free non-elementary Gromov-hyperbolic groups \cite{Har--88}, 
\item"$\circ$"
   $PSL_n(\bold Z)$ for $n \ge 2$ \cite{BCH1--94},   
   and more generally any Zariski-dense subgroup with centre reduced to
   $\{1\}$ in a connected semi-simple real Lie group 
   without compact factor \cite{BCH2--94},
\item"$\circ$"
   centerless mapping class groups and outer automorphism groups 
   of free groups \cite{BrHa--04},
\item"$\circ$"
   irreducible Coxeter groups which are neither finite nor affine
   (see \cite{Fen--03} and our Corollary~18),
\endroster
as we will discuss below.

\medskip

   Observe that a group $G$ not reduced to one element cannot be both
amenable and C$^*$-simple. Indeed, this would imply $\lambda_{G} \prec 1_{G}$; 
but any coefficient of $1_{G}$ is a constant function on $G$,  
whereas a diagonal matrix coefficient 
$\varphi_{\xi,\xi} : g \longmapsto \langle \xi \mid \lambda_{G}(g) \xi \rangle$,
for an appropriate function $\xi \in L^2(G)$ 
supported in a small neighbourhood $\Cal V$ of the identity,  
has value~$1$ at $g = 1$ and~$0$ at
$g$  in the complement of $\Cal V$; 
thus $\varphi_{\xi,\xi}$ is not a limit of constant functions, 
so that $\lambda_G \not\prec 1_G$. \par

   Together with a basic result on unitary induction of representations, 
the same argument is used to show the following generalization.

\medskip

\proclaim{3.~Proposition} 
Let $G$ be a second-countable 
\footnote{
To avoid technical problems, we will assume whenever useful
that locally compact groups appearing here are {\it second countable.} 
Accordingly, \lq\lq discrete groups\rq\rq , 
denoted by $\Gamma$, are assumed to be {\it countable}.
}
locally compact group
which contains an amenable closed normal subgroup $N \neq \{1\}$. 
Then $G$ is not C$^*$-simple. 
\endproclaim

\demo{Proof} We formulate the proof 
for a group $\Gamma$ which is discrete and countable.
The general case is similar, 
modulo a standard (but technical) result 
on the existence of quasi-invariant measures 
on  homogeneous spaces of $G$. \par

   Consider the {\it quasi-regular representation} 
$\lambda_{\Gamma/N}$ of $\Gamma$ in $\ell ^2 (\Gamma / N)$,
defined by 
$$
\left(\lambda_{\Gamma/N}(\gamma)\xi \right) (\gamma 'N) =
\xi(\gamma^{-1} \gamma' N)
\qquad \text{for all} \qquad
\gamma \in \Gamma, 
\quad \xi \in \ell^2(\Gamma/N),
\quad \gamma ' N \in \Gamma /N ;
$$
it coincides with the so-called
{\it induced repre\-sentation} 
$\operatorname{Ind}_N^{\Gamma}(1_N)$ of $1_N$ from $N$ to $\Gamma$.
Observe that any diagonal coefficient of $\lambda_{\Gamma/N}$ is constant on $N$; 
indeed, for $\xi \in  \ell ^2(\Gamma / N)$ 
and  $\gamma \in N$,  we have
$\big( \lambda_{\Gamma/N}(\gamma) \xi \big) (xN) = \xi(\gamma^{-1}xN) = \xi (xN)$
for all $x \in \Gamma$, namely $\lambda_{\Gamma/N}(\gamma)\xi = \xi$.
\par

   On the contrary, diagonal matrix coefficients of $\lambda_{\Gamma}$
separate the points of $\Gamma$.
Indeed, let $a \in \Gamma$; set $\xi = 2^{-1/2}(\delta_1 + i \delta_a)$;
then $\langle \xi \mid \lambda_{\Gamma}(a) \xi \rangle = 1$ when $a = 1$
and $\langle \xi \mid \lambda_{\Gamma}(a) \xi \rangle \ne 1$ when $a \ne 1$.

   We have $1_N \prec \lambda_N$ by hypothesis. 
As weak containment is stable by induction (a result of Fell \cite{Fel--64}),  
we have also 
$$
\operatorname{Ind}_N^{\Gamma}(1_N) = \lambda_{\Gamma/N} 
\ \prec \ 
\lambda_{\Gamma} = \operatorname{Ind}_N^{\Gamma}(\lambda_N) .
$$
If $\Gamma$ was C$^*$-simple, 
we would have $\lambda_{\Gamma} \prec \lambda_{\Gamma/N}$; 
but this is impossible since we have just checked that 
any diagonal coefficient of $\lambda_{\Gamma/N}$ is constant on $N$, 
whereas this does not  hold for $\lambda_{\Gamma}$. 
\hfill $\square$ 
\enddemo

\demo{Reformulation}
In any locally compact group $G$, 
there exists a unique maximal amenable normal subgroup 
which is called the {\it amenable radical}.
The notion is due to Day (see Lemma~1 in \S~4 of \cite{Day--57}),
and a good reference is Proposition~4.1.12 in \cite{Zim--84}.
Thus Proposition~3 can be reformulated as follows:
{\it
If $G$ is a second-countable locally compact group 
with amenable radical not reduced to one element, 
then $G$ is not C$^*$-simple.
}
\enddemo

\medskip

\proclaim{4.~Question}
Does there exist a countable group $\Gamma$ 
with amenable radical reduced to one element
such that $\Gamma$ is not C$^*$-simple?
\endproclaim

Question~4 is discussed in \cite{BeHa--00}.

\medskip

 A non-trivial connected group is never C$^*$-simple; 
see Proposition 4  of \cite{BCH1--94}, 
which is a consequence of Proposition 3 above and 
of the structure of connected locally compact groups. 
This explains why, in this exposition, 
\newline
\centerline{
{\it
all positive results on C$^*$-simple groups are about countable groups.
}
}
It is however natural to formulate the following open problem.

\medskip

\proclaim{5.~Question} 
Do there exist non-discrete second countable locally compact 
groups  which are C$^*$-simple?
\endproclaim

\medskip

With Proposition 3, 
representation theory shows two opposite classes of locally compact groups,
namely amenable groups and C$^*$-simple groups.
Let us recall here that the point of view of representations
provides another and even more fundamental tame-versus-chaotic dichotomy
between so-called {\it groups of type I} and other groups
(it is important here that
locally compact groups are {\it second-countable}).
See Part VI of the appendix for the relevant theorem, essentially due Glimm, 
and Part VII for groups of \lq\lq type I\rq\rq .

If $G$ is of type I, 
the set of its irreducible unitary representations
has a Borel structure which is countably separated, 
and its representations can be analyzed 
as \lq\lq direct integrals\rq\rq \ of irreducible representations
in a useful way.
Two irreducible representations of such a $G$  
are weakly equivalent if and only they are unitarily equivalent.
Groups of type I include 
compact groups,
abelian groups, 
connected real Lie groups which are either semi-simple
or nilpotent, 
and many others (more on this in the Appendix).

On the contrary, for {\it any group which is not of type I,} 
there is no reasonable description 
\footnote{
The statement on \lq\lq no reasonable description\rq\rq \ may be unfair
in the case of nilpotent groups; 
see for instance \cite{Kir--04} and \cite{BaKM--97}.
}
of all irreducible unitary representation;
moreover, there is no kind of unicity in the decomposition of a unitary
representations in irreducible representations.
Groups not of type I include all countable groups 
which are not virtually abelian 
and some soluble connected real Lie groups 
(there are examples in \S~19 of \cite{Kir--74}).

\bigskip
\head
{\bf
2. Formulation in terms of C$^*$-algebras
} 
\endhead

   Let $\Gamma$ be a group. To any unitary representation 
$\rho : \Gamma \longrightarrow \Cal U (\Cal H)$, we
associate the C$^*$-algebra
$$
   C^*_{\rho}(\Gamma) \, = \, 
   \text{norm closure of linear span of } \
   \big(  \rho(\gamma) \big) _{\gamma \in \ \Gamma}
$$
which is a sub-C$^*$-algebra of the algebra of all bounded operators on $\Cal H$.
More generally, for a unitary representation $\rho$ of a second countable locally
compact group $G$, we define 
$$
   C^*_{\rho}(G)  \, = \, \text{norm closure of linear span of } \ 
       \int_G \rho(g) f(g) dg \ \text{ for } \ f \in \Cal K (G)
$$
where $\Cal K (G)$ denotes the convolution algebra of complex-valued continuous
functions with compact supports on $G$, and where  $dg$ denotes a
left-invariant Haar measure on $G$. \par

  In the trivial case where $\rho$ is one-dimensional, we find
$C^*_{\rho}(G) \approx \bold C$. In case $\rho$ is irreducible and finite
dimensional, say of dimension $n$, we have 
$C^*_{\rho}(G) \approx Mat_n(\bold C)$ by a theorem of Burnside.
(See e.g. \S \ 6 in \cite{Ser--67}, or \S \ XVII.3 in \cite{Lan--65}, 
at least for the case of a finite group.) 
Besides these \lq\lq small\rq\rq \ examples, 
the two most important cases appear in the next definition.
Recall that the  {\it universal representation} $\pi_{uni}$ of a group $G$
is the direct sum of all (equivalence classes of) 
unitary representations of $G$ in appropriate Hilbert spaces
\footnote{
Some bound should be given on the dimensions of these spaces,
so that the sum is over a well-defined set.
}
\hskip-.1cm
.

\medskip

\proclaim{6.~Definition}
The {\rm maximal C$^*$-algebra} $C^*_{max}(G)$ of $G$
is the C$^*$-algebra corresponding 
to the  {\rm universal representation} $\pi_{uni}$ of $G$. 
The {\rm reduced C$^*$-algebra}  $C^*_{red}(G)$ of $G$
is the C$^*$-algebra corresponding 
to the left-regular representation of $G$.
\endproclaim

\medskip

The maximal C$^*$-algebra has some properties of the 
\lq\lq complex group algebra\rq\rq \ in algebra
\footnote{
In the case of \lq\lq discrete groups\rq\rq , it is also functorial in the
following sense : any group homomorphism $\Gamma_1 \longrightarrow \Gamma_2$ 
\lq\lq extends\rq\rq \ naturally to a morphism of C$^*$-algebras 
$C^*_{max}(\Gamma_1) \longrightarrow C^*_{max}(\Gamma_2)$,
and the C$^*$-morphism is injective if the group morphism is injective
(Proposition 1.2 of \cite{Rie--74}).
This does {\it not} carry over to locally compact groups !
For example, $C^*_{max}(PSL_2(\bold Z))$, a C$^*$-algebra with unit,
does not embedd in $C^*_{max}(PSL_2(\bold R))$, which has none.
For a lattice $\Gamma$ in a Lie group $G$, the natural morphism from
$C^*_{max}(\Gamma)$ has for its range the {\it multiplier algebra}
$M(C^*_{max}(G))$; but, whatever functor could be defined here,
it will be such that this morphism  
$C^*_{max}(\Gamma) \longrightarrow M(C^*_{max}(G))$
needs not be injective \cite{BeVa--95}.
}
\hskip-.05truecm
. 
In particular,
any unitary representation
$\rho : G \longrightarrow \Cal U (\Cal H_{\rho})$  
determines naturally a morphism of C$^*$-algebras  
$C^*_{max}(G) \longrightarrow C^*_{\rho}(G)$,
again denoted by $\rho$, with kernel denoted by $C^*Ker(\rho)$.
We will see below that 
the reduced C$^*$-algebra is very far from being functorial,
since many residually finite groups
\footnote{
A group $\Gamma$ is {\it residually finite} if, 
for every $\gamma \in \Gamma$, $\gamma \ne 1$,  
there exists a finite group $Q$ and a homomorphism 
$\phi : \Gamma \longrightarrow Q$ 
such that $\phi(\gamma) \ne 1$; 
or, equivalently, if the intersection 
of the normal subgroups of finite index in $\Gamma$ 
is reduced to $\{1\}$. }
have simple reduced C$^*$-algebras. 

\medskip

    The following is standard. See Theorem 1.3 in \cite{Fel--60}, 
or Theorem 3.4.4 and Proposition 18.1.4 in \cite{DC$^*$--69}. 

\medskip

\proclaim{7.~Theorem} Let $G$ be a second countable locally compact group and
let $\rho,\sigma$ be two unitary representations of $G$. 
Then the following properties are equivalent. \smallskip

   $(i)$ \phantom{iii} $C^*Ker (\sigma) \subset C^*Ker (\rho)$, \par
   $(ii)$ \phantom{ii} $\norm{\rho(x)} \le \norm{\sigma(x)}$ for all
              $x \in C^*_{max}(G)$, \par
   $(iii)$ \phantom{i} $\rho \prec \sigma$.
\endproclaim

\demo{On the proof} 
If $(i)$ holds, the morphism of C$^*$-algebras
$$
\left\{ \aligned
C^*_{\sigma}(G) \, = \, C^*_{max}(G) / C^*Ker(\sigma) \
   &\longrightarrow \ 
   C^*_{max}(G) / C^*Ker(\rho) \, = \, C^*_{\rho}(G) \\
\sigma(x) \ &\longmapsto \ \rho(x)
\endaligned \right.
$$
is well defined. As morphisms of C$^*$-algebras are contractions, $(ii)$ follows.
Conversely, it is obvious that $(ii)$ implies $(i)$. 

\medskip

   The proof of the equivalence of $(i)$ and $(ii)$ with $(iii)$ requires some
prerequisite. Recall that a function $\varphi : G \longrightarrow \bold C$ is of 
{\it positive type} if
$\sum_{1 \le j,k \le n} \varphi \left( g_j^{-1}g_k \right) 
\overline{\lambda_j} \lambda_k \ge 0$
for all $g_1,\hdots,g_n \in G$ and $\lambda_1,\hdots,\lambda_n \in \bold C$.
If $\rho$ is a unitary representation of $G$, it is straightforward that,
for any vector $\xi$ in the representation space of $\rho$, the function
$\varphi^{\rho}_{\xi} : g \longmapsto \langle \xi \mid \rho(g) \xi \rangle$
is of positive type. It is a fact that any function of positive type $\varphi$ on
$G$ extends to a positive linear form on $C^*_{max}(G)$; 
if $\varphi = \varphi^{\rho}_{\xi}$ is associated to $\rho$, this linear form
vanishes on $C^*Ker(\rho)$, and defines consequently a positive linear form on
$C^*_{\rho}(G)$. Moreover, a sequence $\left( \varphi_i \right)_{i \ge 1}$ 
of functions of positive type on $G$ converges uniformly on
compact subsets of $G$ if and only if the corresponding sequence of bounded 
linear forms converges in the weak-$*$ topology on 
the appropriate C$^*$-algebra. 
For all this, see Sections 2.7, 13.4, 13.5, and 18.1 in \cite{DC$^*$--69}. 

\medskip

   We can now sketch a proof that the negation of $(iii)$ implies the negation
of $(i)$. If $\rho \nprec \sigma$, there exist some positive linear form
$\varphi^{\rho}_{\xi} : C^*_{\max}(G) \longrightarrow \bold C$
which is {\it not} in the weak-$*$-closure of the convex hull of the positive
linear forms $\varphi^{\sigma}_{\eta}$ associated to $\sigma$.
It is then a consequence of the Hahn-Banach theorem that there exists
$x \in C^*_{\max}(G)$ such that 
$\langle \eta \mid \sigma(x) \eta \rangle = 0$
for all $\eta \in \Cal H_{\sigma}$, and 
$\langle \xi \mid \rho(x) \xi \rangle \neq 0$
for some $\xi \in \Cal H_{\rho}$. 
Hence
\footnote{
Recall that an operator $y$ on a {\it complex} Hilbert space $\Cal H$ is zero if
and only  if $\langle \xi \mid y \xi \rangle = 0$ for all $\xi \in \Cal H$;
see e.g. \S \ 18 in \cite{Hal--51}.
This does {\it not} hold for real Hilbert spaces, as the matrix
$\left( \matrix 0 & 1 \\ -1 & 0 \endmatrix \right)$ demonstrates.
}
$\sigma(x) = 0$ and $\rho(x) \neq 0$, so that the negation of $(i)$ holds.  

\medskip

   To show that $(iii)$ implies $(ii)$, one can argue as follows;
we assume for simplicity that $G = \Gamma$ is a countable group, 
namely that the C$^*$-algebras which occur have  units.
For $\xi \in \Cal H_{\rho}$, the positive linear form
$\varphi^{\rho}_{\xi}$ on $C^*_{max}(\Gamma)$ 
has a norm which is its value at $1$ :
$$ 
\norm{ \varphi^{\rho}_{\xi} } \, = \, \langle \xi \mid \rho(1) \xi \rangle
   \, = \, \norm{\xi}^2 .
$$
Suppose now that $\rho \prec \sigma$. For each compact subset
$Q$ of $G$ and for each $\epsilon > 0$, there exist a finite sequence
$\eta_1,\hdots,\eta_n$ of vectors in $\Cal H_{\sigma}$ and a corresponding 
positive linear form
$$
   \psi^{\sigma}_{Q,\epsilon} \, : \,   x \, \longmapsto \, \sum_{j=1}^N 
              \langle \eta_j \mid \sigma(x) \eta_j \rangle
$$
on $C^*_{max}(\Gamma)$, 
of norm $\psi^{\sigma}_{Q,\epsilon}(1) = \sum_{j=1}^N\norm{\eta_j}^2$, 
such that the net $\left( \psi^{\sigma}_{Q,\epsilon} \right)_{Q,\epsilon}$ 
converges to $\varphi^{\rho}_{\xi}$ in the weak*-topology.
[The set of pairs $(Q,\epsilon)$ with $Q$ compact in $G$ and $\epsilon > 0$
is a net with respect to the order
$(Q,\epsilon) \le (Q',\epsilon ')$ if $Q \subset Q'$ and $\epsilon \ge \epsilon '$.]
Thus,
for $(Q,\epsilon)$ in the net and 
$\psi^{\sigma}_{Q,\epsilon}(^.) = 
\sum_{j=1}^N \langle \eta_j \mid \sigma(^.) \eta_j \rangle$,
we have
$$
\aligned
 \vert \langle \xi \mid \rho(x) \xi \rangle \vert 
 \, &\le \, 
 \sum_{j=1}^N \vert \langle \eta_j \mid \sigma(x) \eta_j \rangle \vert 
    \, + \, \epsilon
 \, \le \, \norm{ \sigma (x) } \sum_{j=1}^N \norm{\eta_j}^2 \, + \, \epsilon
\\
 \, &\le \, 
 \norm{ \sigma (x) } \, \norm{\xi}^2 \, + \, 2 \epsilon
\endaligned
$$
for all $x \in C^*_{max}(\Gamma)$. 
This shows that 
$\vert \langle \xi \mid \rho(x) \xi \rangle \vert  \le 
\norm{ \sigma (x) } \, \norm{\xi}^2$, 
and thus also that
$\norm{ \rho(x) } \le \norm{ \sigma (x) }$, for all $x \in C^*_{max}(\Gamma)$. 
\hfill $\square$
\enddemo

\medskip

As a consequence, we can justify the terminology of Section 1. 

\medskip

\proclaim{8.~Corollary} A second-countable locally compact group $G$ is
C$^*$-simple if and only if its reduced C$^*$-algebra is simple. 
\endproclaim

\demo{Proof} Assume $C^*_{red}(G)$ is simple as an algebra. Then
$C^* Ker(\lambda_G)$ is maximal as a two-sided ideal in $C^*_{max}(G)$, hence
$G$ is C$^*$-simple in the sense of Definition~2. \par

   Assume conversely that $G$ is C$^*$-simple as in Definition~2, 
so that $C^*Ker(\lambda_G)$ is maximal among ideals of the form $C^*Ker(\rho)$.
Then $C^*Ker(\lambda_G)$ is maximal among all two-sided ideals of $C^*_{max}(G)$
by a general result on C$^*$-algebras (Theorem 2.9.7.ii in \cite{DC$^*$--69}),
and consequently $C^*_{red}(G)$ is a simple algebra. 
\hfill $\square$
\enddemo

\medskip

   Note that the {\it maximal} C$^*$-algebra of a group $G \ne 1$ 
cannot be  simple since the unit representation $1_G$
provides a morphism from $C^*_{max}(G)$  onto $\bold C$. 
\par

   Let $G$ be a second-countable locally compact group. 
If $G$ is amenable, {\it any} unitary  representation of $G$ 
is weakly contained  in the regular one
(as already stated after Proposition~2), 
so that the canonical morphism from
$C^*_{max}(G)$ onto $C^*_{red}(G)$ is an isomorphism. 
In the particular case of an abelian group, 
$C^*_{red}(G)$ is also isomorphic to the algebra 
$\Cal C _o (\hat G)$ of continuous function which vanish at infinity 
on the dual $\hat G$ of $G$ (Fourier transform and Pontrjagin duality); 
it follows that closed ideals of $C^*_{red}(G)$ 
are in a natural bijective correspondance with closed subspaces of $\hat G$. 
\par

   It is a natural project to investigate in general 
the structure of closed ideals 
in $C^*_{max}(G)$ and $ C^*_{red}(G)$. 
Let us indicate one motivation.
Let $A$ be a C$^*$-algebra, let $I$ be a closed two-sided ideal,
and let $A/I$ denote the quotient C$^*$-algebra.
Some problems concerning $A$ split as simpler problems for $I$ and $A/I$.
For example, this is sometimes the case for the computation of K-groups,
because there exists a cyclic six-term exact sequence connecting $K_0$ and $K_1$
of $I$, $A$, and $A/I$.
Recall that the K-groups of algebras of the form $C^*_{red}(\Gamma)$
are crucial ingredients of the {\it Baum-Connes conjecture}
\cite{Val--02}.
Thus a non-trivial closed two-sided ideal in $A = C^*_{red}(\Gamma)$
can suggests a strategy of computation of $K_*(A)$;
in other words, from this point of view, 
the C$^*$-simple groups provide  \lq\lq the hardest case\rq\rq  .

\medskip

   Observe that Theorem 7 implies various numerical (in)equalities.
For example, if $S$ is a finite subset, say of size $n$, 
of a locally compact group $G$, and if $\rho$ is 
a unitary representation of $G$ which is weakly equivalent to $\lambda_G$, 
then
$$
  \norm{ \frac{1}{n} \sum_{s \in S} \rho(s) } \, = \, 
  \norm{ \frac{1}{n} \sum_{s \in S} \lambda_G(s) } .
$$
See Item 20  below.

\bigskip
\head
{\bf
3. Powers result on non-abelian free groups
}
\endhead

In 1975, Powers established that non-abelian free groups are C$^*$-simple.
Once conveniently reformulated, 
his method proved to be robust enough to apply to many more groups.
The following terminology refers to the original paper \cite{Pow--75}.

\proclaim{9.~Definition} 
A group $\Gamma$ has the {\rm Powers property}, or is a {\rm Powers groups},
if \smallskip
   for any finite subset $F$ in $\Gamma \setminus \{1\}$ and for any integer 
          $N \ge 1$, \par
   there exists a partition $\Gamma = C \sqcup D$ and elements
          $\gamma_1,\hdots,\gamma_N$ in $\Gamma$ such that \par
   $fC \cap C = \emptyset$ for all $f \in F$
          and $\gamma_jD \cap \gamma_kD = \emptyset$ 
          for all $j,k \in \{1,\hdots,N\}$, $j \neq k$. 
\endproclaim

\medskip

\proclaim{10.~Definitions} 
A homeomorphism $\gamma$ of a Hausdorff space $L$ is
{\rm hyperbolic} if it has two fixed points $s_{\gamma},r_{\gamma} \in L$ 
such that the following holds : for any neighbourhood $S$ of $s_{\gamma}$  
and $R$ of $r_{\gamma}$, there exists $n_0 \in \bold N$ such that
$\gamma^n(L \setminus S) \subset R$ and 
$\gamma^{-n}(L \setminus R) \subset S$ for all $n \ge n_0$. 
The points $s_{\gamma}$ and $r_{\gamma}$ are  called the
{\rm source} and the {\rm  range} of $\gamma$, respectively.
\par

   Two hyperbolic homeomorphisms of $L$ are {\rm transverse}
if they have no common fixed point. 
\par

   An action of a group $\Gamma$ on $L$ is {\rm strongly faithful} if, 
for any finite subset $F$ of  $\Gamma \setminus \{1\}$, 
there exists $y \in L$ such that $fy \ne y$ for all $f \in F$.
It is {\rm strongly hyperbolic} if, 
for any integer $N \ge 1$, there exist $N$ pairwise transverse hyperbolic 
elements in $\Gamma$. 
\par
\endproclaim

    In the definition of \lq\lq strongly hyperbolic\rq\rq , observe that the
condition for $N = 2$ implies the condition for any $N \ge 2$.
Indeed, if $h_1,h_2$ are pairwise transverse hyperbolic homeomorphisms, 
then $h_1^{j_k} h_2 h_1^{-j_k}$ are pairwise transverse and hyperbolic 
for an appropriate sequence of integers $\left(j_k\right)_{k \ge 1}$
with $\lim_{k \to \infty} j_k = \infty$.

\medskip

\proclaim{11.~Proposition} Let $\Gamma$ be a group acting by homeomorphisms on
a Hausdorff space $L$. Assume the action is 
\roster
\item"(i)" minimal, \par
\item"(ii)" strongly faithful, \par
\item"(iii)" strongly hyperbolic. 
\endroster
Then $\Gamma$ is a Powers group.
\endproclaim

\demo{Proof} Consider $F \subset \Gamma \setminus \{1\}$ and $N \ge 1$ as in the
definition of \lq\lq Powers group\rq\rq . By $(ii)$ there exists a non-empty open
subset $C_L$ in $L$ such that $f C_L \cap C_L = \emptyset$ for all $f \in F$.
By $(iii)$ there exist $N+1$ pairwise transverse hyperbolic elements 
$\gamma = \gamma_0,\gamma_1,\hdots,\gamma_n$ in $\Gamma$. By $(i)$ we may assume
$r_{\gamma} \in C_L$. \par

   Upon conjugating $\gamma_1,\hdots,\gamma_N$ by a large enough power of
$\gamma$, we may assume that, for each $j \in \{1,\hdots,N\}$, both the source
$s_j$ and the range $r_j$ of $\gamma_j$ are in $C_L$. We choose neighbourhoods
$S_j$ of $s_j$ and $R_j$ of $r_j$ in such a way that
$S_1, R_1, \hdots, S_N, R_N$ are pairwise disjoint and inside $C_L$.
Upon replacing each of $\gamma_1,\hdots,\gamma_N$ by a large enough power of
itself, we may furthermore assume that $\gamma_j(L \setminus C_L) \subset R_j$,
and in particular that $\gamma_j(L \setminus C_L)$ are pairwise disjoint subsets 
of $L$. \par

   Choose $y \in L$ and let 
$C = \left\{ \gamma \in \Gamma \mid \gamma y \in C_L \right\}$ and 
$D = \left\{ \gamma \in \Gamma \mid \gamma y \notin C_L \right\}$. 
Then $fC \cap C = \emptyset$,  since
$fC_L \cap C_L = \emptyset$,  
for all $f \in F$, and $\gamma_j D \cap \gamma_k D = \emptyset$, 
since $R_j \cap R_k = \emptyset$, for all $j \neq k$. 
\hfill $\square$
\enddemo

\medskip

   {\it Note.} In \cite{BrHa--04}, 
Conditions (i) to (iii) have been changed as follows,
so that minimality needs not to be assumed. 
Let $L_{\operatorname{hyp}}$ denote the subset of $L$  
of those points which are fixed by some hyperbolic element in $\Gamma$;
this is a $\Gamma$-invariant subset of $L$. It is enough to assume:
\roster
\item"(i)"
  $\Gamma$ contains two transverse hyperbolic homeomorphisms of $L$;
\item"(ii)"
  For any finite subset $F$ of $\Gamma \setminus \{1\}$,
  there exists a point $t \in L_{\operatorname{hyp}}$
  such that $ft \ne t$ for all $f \in F$.
\endroster

\medskip

\proclaim{12.~Corollary} The following are Powers groups: 
\roster
\item"(i)" free products $\Gamma_1 \ast \Gamma_2$
  with $(\vert \Gamma_1 \vert - 1)(\vert \Gamma_2 \vert - 1) \ge 2$;
\item"(ii)"   non-soluble subgroups of $PSL_2(\bold R)$;
\item"(iii)"  torsion-free non-elementary Gromov-hyperbolic groups;
\item"(iv)" sufficiently large subgroups in mapping class groups of surfaces
of genus $g \ge 1$;
\item"(v)" appropriate subgroups of the group of outer automorphisms of non
abelian free groups.
\endroster
\endproclaim

\demo{Observation} Non-abelian free groups are particular cases of each of
the five classes in the corollary. 
\enddemo

\demo{Proof} For a group in (i), the corollary follows from the proposition
applied to the action of the group on the boundary of the tree
associated to the free product as in \cite{Ser--83}.
For a group in (ii) [respectively in (iii)],
the same argument applies to the action 
on its limit set in the boundary of the
hyperbolic plane on which $PSL_2(\bold R)$ acts by isometries
[respectively on  the Gromov boundary of the group]. 

   For (iv), consider the mapping class group $\tilde \Gamma_g$ of an
orientable closed surface of genus $g \ge 1$, viewed as acting on the
boundary $P\Cal M \Cal F_g$ of the corresponding Teichm\"uller space.
For $g \le 2$, the group $\tilde \Gamma_g$ has a centre $Z_g$ of order $2$
which acts on $P\Cal M \Cal F_g$ as the identity. 
We set $\Gamma_g = \tilde \Gamma_g/Z_g$ for $g \le 2$ 
and $\Gamma_g = \tilde \Gamma_g$ for $g \ge 3$, 
and we consider from now on the natural action 
of $\Gamma_g$ on $P\Cal M \Cal F_g$. (If $g = 1$, there are standard
identifications of $\Gamma_1$ with $PSL_2(\bold Z)$ and of
$P\Cal M \Cal F_1$ with the boundary of the hyperbolic plane.)
It is well-known that the action of $\Gamma_g$ on $P\Cal M \Cal F_g$
is minimal for all $g \ge 1$; see e.g. \S \ VII of expos\'e 6 in
\cite{FLP--79}.

   A {\it pseudo-Anosov class} in $\Gamma_g$ is an element which defines a
hyperbolic homeomorphism of $P\Cal M \Cal F_g$. It is known that $\Gamma_g$
contains transverse pairs of pseudo-Anosov classes; this is for example an
immediate consequence of Lemma 2.5 in \cite{McP--89}. 
Let $\Gamma$ be a subgroup of $\Gamma_g$ which is {\it sufficiently large}
(in the sense of \cite{McP--89}),
which means here that is contains a transverse pair of pseudo-Anosov classes. 
Then there is a unique subset $L_{\Gamma}$ in $P\Cal M \Cal F_g$
which is closed, invariant by $\Gamma$, and minimal for these properties; 
it is the {\it limit set} of $\Gamma$.
If $L^0_{\Gamma}$ is the subset of $P\Cal M \Cal F_g$ 
of those points which are fixed by some pseudo-Anosov class in $\Gamma$, 
then $L_{\Gamma}$ is the closure of $L^0_{\Gamma}$ 
(see Theorem 4.1 in \cite{McP--89}). 
For example and as recalled above, we have $L_{\Gamma} = P\Cal M \Cal F_g$ 
if $\Gamma = \Gamma_g$. It follows from the definitions that the action of
the group $\Gamma$ on its limit set $L_{\Gamma}$ is minimal and strongly
hyperbolic. 

   We claim that the action of $\Gamma$ on $L_{\Gamma}$ 
is also strongly faithful. 
Indeed, let $F$ be a finite subset of $\Gamma$, not containing the identity. 
If $m$ is the number of elements of $F$, choose an integer $n > \frac{m}{2}$
and pairwise transverse pseudo-Anosov classes $h_1,\hdots,h_n \in \Gamma$.
Let $X$ denote the $2n$-element subset of $L_{\Gamma}$ consisting of the
points fixed by one of the $h_j$ 's. Because of Lemma 2.5 in \cite{McP--89}
again, any element of $\Gamma_g$ can fix at most one point in $X$.
As $\abs{X} = 2n > m = \abs{F}$, there is at least one point $x_0 \in X$ such
that $fx_0 \ne x_0$ for all $f \in F$. This ends the proof of the claim, and
thus that of (iv). 
This proof of (iv) is that of  \cite{BrHa--04}. 

We refer to the same paper 
for a precise formulation and for a proof of~(v).
\hfill $\square$
\enddemo

\medskip

\proclaim{13.~Questions}
(i) A Powers group is clearly non-amenable, and is icc 
\footnote{
A group $\Gamma$ has {\it infinite conjugacy classes,} or shortly is {\it icc,} 
if it is infinite 
and if all its conjugacy classes distinct from $\{1\}$ are infinite.  
More on icc groups on Part~X of the appendix.
}
(see Proposition~1 in \cite{Har--85}). 
Does there exist Powers groups without  non-abelian free subgroups? 
\par

 (ii) Let $\Gamma$ be a direct product of two non-abelian free group. Does there
exist a faithful action of $\Gamma$ by homeomorphisms on a Hausdorff space such
that at least one element of the group acts hyperbolically? \par

  (iii) Does there exist an action of $SL_3(\bold Z)$ by homeomorphisms on a 
Hausdorff space such that at least one element of the group acts hyperbolically?
\endproclaim  

   Questions (ii) and (iii) are possibly related 
to the notions of {\it Dynkin spaces}  and {\it $\Cal D$-groups}, 
defined by Furstenberg in \cite{Fur--67}; 
but I have not understood enough of this paper to use it.

\medskip

\proclaim{14.~Theorem (Powers)} 
Powers groups are C$^*$-simple.
\endproclaim

\demo{Proof, which follows that of \cite{Pow--75}} Consider a Powers group
$\Gamma$, its left-regular representation $\lambda$, a non-zero ideal $\Cal I$
of its reduced C$^*$-algebra, and an element $U \ne 0$ in $\Cal I$. We want to
show  that $\Cal I$ contains an element $Z$ such that $\norm{Z - 1} < 1$, and in
particular such that $Z$ is invertible (with inverse $\sum_{n=0}^{\infty}Z^n$).\par

   Upon replacing $U$ by a scalar multiple of $U^*U$, we may assume that
$U = 1 + X$ and $X = \sum_{x \in \Gamma, x \neq 1} z_x \lambda(x)$,
with $z_x \in \bold C$. Choose $\epsilon > 0$, with $\epsilon < 1$; 
there exists a finite subset  $F$ of $\Gamma \setminus \{1\}$ such that, if
$$
   X' \, = \, \sum_{f \in F} z_f \lambda (f) ,
$$
then $\norm{X' - X} < \epsilon$. 
Set $U' = 1 + X'$, so that $\norm{ U' - U} < \epsilon$.
Choose an integer $N$ so large that
$\frac{2}{\sqrt N} \norm{X'} < 1 - \epsilon$. 
\par

   Let now $\Gamma = C \sqcup D$ and $\gamma_1,\hdots,\gamma_N$ be as in the
definition of the Powers property. Set
$$
\aligned
 V \, = \,  \frac{1}{N} \sum_{j=1}^N \lambda(\gamma_j) U  \lambda(\gamma_j^{-1})
 \qquad\qquad &
 V' \, = \, \frac{1}{N} \sum_{j=1}^N \lambda(\gamma_j) U' \lambda(\gamma_j^{-1})
 \\
 Y \, = \,  \frac{1}{N} \sum_{j=1}^N \lambda(\gamma_j) X  \lambda(\gamma_j^{-1})
 \qquad\qquad &
 Y' \, = \, \frac{1}{N} \sum_{j=1}^N \lambda(\gamma_j) X' \lambda(\gamma_j^{-1})
\endaligned
$$
so that $V = 1+Y \in \Cal I$ and $V' = 1 + Y'$. We show below that 
$\norm{Y'} < 1 - \epsilon$. This implies that 
$\norm{Y} \le \norm{Y'} + \norm{Y-Y'} \le \norm{Y'} + \norm{X-X'} < 1$.
As $\Cal I$ contains the invertible element $V = 1 + Y$, 
the C$^*$-algebra $C^*_{red}(\Gamma)$ is indeed simple.  \par

   For $j \in \{1,\hdots,N\}$, denote by $P_j$ the orthogonal projection of 
$\ell ^2 (\Gamma)$ onto $\ell ^2 (\gamma_j D)$. 
We have
$$
   (1-P_j) \lambda(\gamma_j) X' \lambda(\gamma_j^{-1}) (1-P_j)
   \, = \, 0 ;
$$
indeed, since $fC \cap C = \emptyset$ for all $f \in F$, we have
$$
\aligned
 \bigg( \lambda(\gamma_j)X'\lambda(\gamma_j^{-1})(1-P_j) \bigg) (\ell^2(\Gamma))
 \ &\subset \
 \bigg( \lambda(\gamma_j)X' \bigg) (\ell^2(C))
\\
 \ &\subset \
 \bigg( \lambda(\gamma_j) \bigg) (\ell^2(D) ) \ = \ P_j (\ell^2(\Gamma)) .
\endaligned
$$
It follows that
$$
  V' \, = \, 1 \, + \, 
    \frac{1}{N} \sum_{j=1}^N P_j \lambda(\gamma_j) X' \lambda(\gamma_j^{-1})
    \, + \, \left(
    \frac{1}{N} \sum_{j=1}^N P_j \lambda(\gamma_j) X' \lambda(\gamma_j^{-1})
              (1-P_j) \right) ^* .
$$
Since the subsets $\gamma_jD$ of $\Gamma$ are pairwise disjoint, the operators
$X'_j \Doteq P_j \lambda(\gamma_j) X' \lambda(\gamma_j^{-1})$
have pairwise orthogonal ranges in $\ell ^2 (\Gamma)$, and we have
$$
   \norm{ \frac{1}{N} \sum_{j=1}^N X'_j } \, \le \,
   \frac{1}{ \sqrt N } \max_{j=1}^n \norm{ X'_j } \, \le \, 
   \frac{1}{ \sqrt N } \norm{ X' } .
$$
Similarly 
$$
   \norm{ \left( \frac{1}{N} \sum_{j=1}^N X'_j (1-P_j) \right) ^* } \, = \,  
\norm{ \frac{1}{N} \sum_{j=1}^N X'_j (1-P_j) } \, \le \, \frac{1}{\sqrt N}\norm{X'}.
$$
Consequently
$$
   \norm{Y'} \, = \, \norm{ V' - 1 } \, \le \, \frac{2}{ \sqrt N } \norm{X'}
   \, < \, 1 - \epsilon .
$$
As already observed,  this ends the proof. 
\hfill $\square$
\enddemo

   There is in Section 6 of \cite{BeLo--00} a formulation of the previous proof 
in terms of functions of positive type. 

\medskip

   There are several weak forms of the Powers property 
which imply C$^*$-simplicity. 
One of them has been introduced by Boca and Nitica \cite{BoNi--88},
and has been used to establish that the groups $P_k/C_k$ and $B_k/C_k$
are simple for all $k \ge 3$ (see \cite{GiHa--91} and \cite{BeHa--00});
here, $B_k$ stands for the Artin braid group on $k$ strings,
$P_k$ for its subgroup of pure braids,
and $C_k$ for their common center (which is infinite cyclic).
Other notions, reminiscent of the Power property, 
have been introduced and studied in 
\cite{B\'ed--91}, \cite{B\'ed--93}, \cite{B\'ed--96}, and \cite{BoNi--90}.

\medskip

   Let us end this section by some historical information 
that we have learned from R.~Kadison (see also \cite{Val--89}). 
In a conversation with him of 1949, I. Kaplansky asked whether
any simple C$^*$-algebra with unit other than $\bold C$
has a projection distinct from $0$ and $1$.
In 1968, Kadison suggested to R. Powers to study from this point of view
the reduced C$^*$-algebra $C^*_{red}(F_2)$ of the free group of rank $2$.
Powers showed within a week that it is simple 
(even if he published this seven years later).
But it was only in 1982 that M.~Pimsner and D.~Voiculescu showed that
$C^*_{red}(F_2)$ does not have any non-trivial idempotent \cite{PiV--82}.
Shortly before, Blackadar had shown the first examples of C$^*$-algebras
without non-trivial idempotent \cite{Bla--80, Bla--81}.

\bigskip
\head
{\bf
4. Other examples of C$^*$-simple groups
}
\endhead

\medskip

   Let $G$ be a non-compact simple connected real Lie group and let $\Gamma$ be a
Zariski-dense subgroup of $G$ with centre reduced to one element. If $G$ has real
rank one, we can consider the action of $\Gamma$ on the sphere at infinity of
the hyperbolic space defined by $G$; then $\Gamma$ contains enough  hyperbolic
homeomorphisms, and the argument of Proposition 11 shows that $\Gamma$ is a Powers
group.  
\par

   If $G$ has higher real rank, the action of $\Gamma$ on an appropriate
\lq\lq boundary\rq\rq \ does not give rise to hyperbolic homeomorphisms with
two fixed points each; but there are elements acting in a rather simple way,  
with finitely many fixed points which are sources, sinks and saddle-type points. 
This simple dynamics can be used  to show the following combinatorial property~:
\smallskip

   for any finite subset $F$ of $\Gamma \setminus \{1\}$, \par
   there exist $y_0 \in \Gamma$ and subsets $U,A_1,\hdots,A_N$ in $\Gamma$
such that \par
   (i) \phantom{iii} $\Gamma \setminus U \subset A_1 \cup \hdots \cup A_N$, \par
   (ii) \phantom{ii} $fU \cap U = \emptyset$ for all $f \in F$, \par
   (iii) \phantom{i} $y_0^{-j}A_s \cap A_s = \emptyset$ for all $j \ge 0$
       and $s \in \{1,\hdots,N\}$.

\smallskip\noindent And this property is sufficient to imply
the following result, from \cite{BCH2--94}.

\medskip

   \proclaim{15.~Theorem} Let $G$ be a connected real semisimple Lie group without
compact factors. Let $H$ be a subgroup of $G$ with trivial centre, whose image in
the adjoint group of the Lie algebra of $G$ is Zariski-dense, and let $\Gamma$ be
the group $H$ together with the discrete topology. Then $\Gamma$ is C$^*$-simple.
\endproclaim

   Theorem 15 shows for example that $PSL_n(\bold K)$ is C$^*$-simple for any
$n \ge 2$ and for $\bold K = \bold Q$, or more generally for $\bold K$ a number
field.

\medskip

   The proof of Theorem 15 which is presently available is substancially longer
than  that of Powers Theorem 14. A positive answer to the following question would
provide a shorter proof, and be of independent interest. (For simplicity, we
consider here torsion-free groups; for groups which have torsion, see
\cite{BCH2--94}.) 

\medskip

\proclaim{16.~Question} 
Let $\Gamma$ be a group as in Theorem 15. Is it true that,
for any finite subset $F$ of $\Gamma \setminus \{1\}$, there exists $y \in \Gamma$
such that the subgroup $\langle x,y \rangle$ of $\Gamma$ generated by $x$ and $y$ is free
of rank two for all $x \in F$ ? 
\endproclaim

\noindent The answer is \lq\lq yes\rq\rq \ if $G$ has real rank one. 
The question makes sense for $G$ a semisimple algebraic group over the $p$-adics.
The next question is related to Question 13.i.

\medskip

\proclaim{17.~Question} 
Does there exist a  countable group which is
C$^*$-simple and which does not contain non-abelian free subgroups? \par
\endproclaim

\medskip

It is essentially a theorem of Fendler,
and also a consequence of Theorem 15,  
that an irreducible Coxeter group which is neither finite nor affine
is C$^*$-simple. 
As both the original paper \cite{Fen--03} and a later paper \cite{BtHa--05}
prove this with a non-necessary extra hypothesis,
let us revisit this result.
Consider
\roster
\item"$\circ$" 
   an irreducible Coxeter system $(W,S)$, with a finite set $S$ of generators,
\item"$\circ$"  
   the free real vector space $E = R^{S}$ on $S$, 
\item"$\circ$" 
   the Tits form $B$ on $E$,
\item"$\circ$" 
   the kernel $E_0 = \{v \in E \mid B(v,w) = 0 \quad \text{for all}
   \quad w \in E\}$ of $B$,
\item"$\circ$" 
   the orthogonal group $Of(E,B)$ of those linear automorphisms $g$ of $E$
   such that $B(gv,gw) = B(v,w)$ for all $v,w \in E$ and
   $gv = v$ for all $v \in E_0$,
\item"$\circ$" 
   the geometric representation $\sigma : W \longrightarrow Of(E,B)$,
   which is faithful (Tits theorem),
\item"$\circ$" 
   the non-degenerate symmetric bilinear form $\tilde B$
   induced by $B$ on the quotient space $\tilde E = E/E_0$,
\item"$\circ$" 
   the corresponding orthogonal group $O(\tilde E, \tilde B)$
   and the canonical projection 
   $\pi$ from $Of(E,B)$ onto $O(\tilde E, \tilde B)$,
\item"$\circ$"
   and the composition $\tilde \sigma : W \longrightarrow O(\tilde E, \tilde B)$
   of $\sigma$ and $\pi$.
\endroster
It is standard that $W$ is finite if and only if $B$ is positive definite.
It is equally standard that $W$ is infinite 
and contains a free abelian group of finite index 
if and only if $B$ is positive and $\dim_{\bold R}(E_0) \ge 1$; 
in this case, $W$ is said to be {\it affine}.
For all this, see for example \cite{Bou--68}.

   Assume moreover that $W$ is neither finite nor affine. 
Fendler has shown that $W$ is C$^*$-simple under the extra condition
that $B$ is non degenerate, namely that $E_0$ is reduced to zero.
The extra assumption is redundant:

\proclaim{18.~Corollary} An irreducible Coxeter group 
which is neither finite nor affine is C$^*$-simple.
\endproclaim

\demo{Proof}
Let the notation be as above.
It is a result of Vinberg that 
the representation $\tilde \sigma$ of $W$ is faithful;
see Proposition 13 in  \cite{Vin--71},
as well as Proposition 6.1.3 in \cite{Kra}, and \cite{Cor}.
It follows from \cite{BtHa--05} 
that its image is Zariski-dense in $O(\tilde E, \tilde B)$.
It is then a consequence of Theorem 15 that $W$ is C$^*$-simple;
details are as in \cite{BtHa--05}.
\hfill $\square$ 
\enddemo

\medskip

Let us  state some hereditary properties of C$^*$-simplicity.

\proclaim{19.~Proposition} $(i)$ The direct product of two C$^*$-simple groups is 
C$^*$-simple. \par
  $(ii)$ An inductive limit of C$^*$-simple groups is C$^*$-simple. \par
   $(iii)$ In a C$^*$-simple group, any subgroup of finite index is C$^*$-simple.
\par
   $(iv)$ Let $\Gamma$ be a group which has a C$^*$-simple subgroup of finite
index. Then $\Gamma$ is C$^*$-simple if and only if it is icc.
 \endproclaim 

\demo{Proof} For $(i)$, see \cite{Tak--64}. For $(ii)$, which is an exercise, 
$(iii)$, and $(iv)$, see \cite{BeHa--00}. 
For $(iii)$ and $(iv)$, see also \cite{Pop--00}. 
\hfill $\square$ 
\enddemo

\medskip

   In all examples of a countable group $\Gamma$ for which it is known that
$C^*_{red}(\Gamma)$ is simple, it is also known that the canonical trace on 
$C^*_{red}(\Gamma)$ is the unique tracial state. Indeed, it is an open problem to 
know whether there exists a group such that $C^*_{red}(\Gamma)$ is simple with
several tracial states, or a group such that $C^*_{red}(\Gamma)$ has a unique
tracial state and is not simple.

   For $\Gamma$ a group as in one of Theorems 14 and 15, 
it is known that a reduced crossed product $A \rtimes_{red} \Gamma$ is simple 
if and only the only $\Gamma$-invariant ideals of the C$^*$-algebra $A$ 
are $0$ and~$A$. 
Thus, for example, if $\Gamma$ is a torsion-free hyperbolic group 
with Gromov boundary~$\partial \Gamma$,
the crossed product $\Cal C (\partial \Gamma) \rtimes_{red} \Gamma$
defined by the action of $\Gamma$ on its boundary is a simple C$^*$-algebra.
More on these algebras in \cite{Ana--97}.

   Observe that, if $\Gamma$ is a group acting on a C$^*$-algebra $A$,
the condition that the only $\Gamma$-invariant ideals in $A$ are $0$ and~$A$
is always necessary for $A \rtimes_{red} \Gamma$ to be simple.
If $\Gamma$ is not a Powers group, this condition needs not be sufficient:
when $\Gamma = \bold Z$ acts on $\bold C$, the reduced crossed product
is $C^*_{red}(\bold Z) \approx \Cal C(\bold S^1)$, which is far from being simple.

\medskip

\subhead
20.~A computational consequence of C$^*$-simplicity
\endsubhead
The knowledge that a group is C$^*$-simple is not only a qualitative result, 
but can sometimes be used in computations. Consider for example the natural
measure-preserving action of $\tilde \Gamma = GL(2,\bold Z)$ on the $2$-torus 
$\bold T ^2$, and the corresponding measure preserving action of 
$\Gamma = PGL(2,\bold Z)$ on the $2$-sphere $\bold S ^2$, which is the quotient of 
$\bold T ^2$ by the action of the centre of $\tilde \Gamma$, as in \cite{Sch--80}.
\par

   The unitary representation $\tilde \pi$ of $\tilde \Gamma$ in the space 
$L^2_0 (\bold T^2) = \left\{ \xi \in L^2(\bold T^2) \mid 
\int_{\bold T^2} \xi(t) dt = 0 \right\}$ 
is weakly contained in the regular representation of $\tilde \Gamma$; to check
this, observe that $L^2_0 (\bold T^2)$ can be seen as an orthogonal direct sum of
$\ell ^2$-spaces corresponding to the orbits of $\tilde \Gamma$ in 
$\bold Z ^2 \setminus \{0\}$, namely of $\ell ^2$-spaces of the form
$\ell ^2 (\tilde \Gamma / N)$ with $N$ amenable. It follows that the 
representation $\pi$ of $\Gamma$ on the space
$L^2_0 (\bold S^2) = \left\{ \xi \in L^2(\bold S^2) \mid 
\int_{\bold S^2} \xi(x) d\nu(x) = 0 \right\}$ 
is weakly contained in $\lambda_{\Gamma}$ (where $d\nu(x)$ is the measure on 
$\bold S ^2$ which is the direct image of 
the Lebesgue measure $dt$ on $\bold T^2$). 
On the other hand, it follows from Proposition 19 that
$PGL(2,\bold Z)$ is C$^*$-simple
(details in \cite{BeHa--00}). \par

   Hence $\norm{\pi(\mu)} = \norm{ \lambda_{\Gamma}(\mu) }$ for any finite 
measure (positive or not) on $\Gamma$. Computations of norms like 
$\norm{ \lambda_{\Gamma}(\mu) }$, at least for positive measures on $\Gamma$, 
go back to Kesten's theorems on random walks on groups \cite{Kes--59}.

\bigskip
\head
{\bf
Appendix
}
\endhead

\subhead 
I.~On the definition of weak containment (Definition~1)
\endsubhead 

\medskip

Let $G$ be a locally compact group, 
$\rho : G \longrightarrow \Cal U (\Cal H)$ and
$\sigma : G \longrightarrow \Cal U (\Cal K)$ two unitary representations,
$\xi \in \Cal H$, 
and $\varphi_{\xi,\xi}$ the corresponding diagonal matrix coefficient of $\rho$;
for simplicity, we assume that $\norm{\xi} = 1$. 
Then the following are equivalent:
\roster
\item"(i)"
   $\varphi_{\xi,\xi}$ is a limit of sums 
   of diagonal matrix coefficients of $\sigma$,
   uniformly on compact subsets of $G$ (as in Section 1),
\item"(ii)"
   $\varphi_{\xi,\xi}$ is a limit of convex combinations 
   of diagonal matrix coefficients of $\sigma$,
   uniformly on compact subsets of $G$.
\endroster
If $\pi$ is moreover irreducible, these conditions are also equivalent to:
\roster
\item"(iii)"
   $\varphi_{\xi,\xi}$ is a limit 
   of diagonal matrix coefficients of $\sigma$,
   uniformly on compact subsets of~$G$.
\endroster
The equivalence of (i) and (ii) is an easy exercise:
if $\varphi_{\xi,\xi}$ is a limit 
of sums $\sum_{j=1}^n \psi_j$,
uniformly on compact subsets of $G$,
we can assume that these compact subsets contain~$1$,
so that $\sum_{j=1}^n \psi_j(1)$ is near $\varphi_{\xi,\xi}(1) = 1$. 
Then $\varphi_{\xi,\xi}$ is also a limit 
of convex combinations $\sum_{j=1}^n c_j \tilde{\psi}_j$, 
with $c_j = \psi_j(1) / (\psi_1(1) + \cdots + \psi_n(1))$
and $\tilde{\psi}_j = \psi_j / \psi_j(1)$.
The proof of the equivalence with (iii) is more technical,
and we refer to Appendix F of \cite{BeHV}.

\bigskip

\subhead
II.~On weak equivalence not preserving irreducibility
\endsubhead

\medskip

Consider a group $\Gamma$ which is not reduced to one element.
The {\it right-regular representation} $\rho_{\Gamma}$ of $\Gamma$
acts on the space $\ell^2(\Gamma)$; it is defined by
$(\rho_{\Gamma}(\gamma)\xi)(x) = \xi(\gamma x)$
for all $\gamma,x \in \Gamma$ and $\xi \in \ell^2(\Gamma)$.
In the algebra of all operators on $\ell^2(\Gamma)$,
the commutant of $C^*_{red}(\Gamma)$ contains $\rho_{\Gamma}(\Gamma)$;
it follows from Shur's lemma that $\lambda_{\Gamma}$ is reducible
\footnote{
A natural question is to ask for 
explicit and elementary constructions of invariant subspaces.  
If $\Gamma$ has at least one element of infinite order,
there is a nice construction of this kind in \cite{Aub--79}, 
in terms of the commutant of  $\lambda_{\Gamma}(\Gamma)$ and of Fourier series. 
I don't know any explicit construction 
in the case of an infinite torsion group.
}
\hskip-.1cm
.
\par

   Assume moreover that $\Gamma$ is {\it icc.} 
Then the representation $\lambda_{\Gamma}$ is factorial,
by Lemma 5.3.4 of \cite{ROIV}.
Recall that a representation $\pi : G \longrightarrow \Cal U (\Cal H)$
of a locally compact group $G$ 
is {\it factorial} if, for any $\pi(G)$-invariant decomposition 
$\Cal H = \Cal H_1 \oplus \Cal H_2$ with $\Cal H_1 \ne \{0\} \ne \Cal H_2$,
there exists a non-zero $G$-intertwining operator from $\Cal H_1$ to $\Cal H_2$;
equivalently, $\pi$ is factorial
if the von Neumann algebra $W^*_{\pi}(G)$ generated by $C^*_{\pi}(G)$
has a center reduced to the multiples of the identity
(see Corollary 5.2.5 of \cite{DC$^*$--69}).
Now, the kernel of any factorial representation of a separable C$^*$-algebra
is primitive, which means that is is also the kernel of an irreducible
representation of the same algebra
(Corollary 3 of Theorem 2 in \cite{Dix--60}, 
or Lemma 2 in \cite{BeHa--00} and Proposition 4.3.6 in \cite{Ped--79}).
As a consequence, we have:

\proclaim{21.~Proposition} For an icc group $\Gamma$, 
the left regular representation $\lambda_{\Gamma}$,
which is reducible,
is weakly equivalent to an irreducible representation.
\endproclaim

Let us provide two specific examples.

\proclaim{22.~Example} Case of a lattice $\Gamma$ in $G = PSL_2(\bold R)$. 
\endproclaim

For $s \in \bold R_+$, let $\pi^{Pr}_s$  denote the corresponding irreducible
representation of the principal unitary series of $G$. It is a particular case
of the results of \cite{CoSt--91} that the restrictions $\pi^{Pr}_s \vert \Gamma$
to the lattice $\Gamma$ are irreducible and pairwise non-equivalent. On the
other hand, $\pi^{Pr}_s \vert \Gamma \prec \lambda_{\Gamma}$, because
$\pi^{Pr}_s \prec \lambda_G$, and weak containment is stable by restriction.
Now for $s \ne s'$, the two representations $\pi^{Pr}_s$, $\pi^{Pr}_{s'}$ are not
weakly equivalent (because they are not unitarily equivalent, and $G$ is
of type I), but $\pi^{Pr}_s \vert \Gamma$ and $\pi^{Pr}_{s'} \vert \Gamma$
are both weakly equivalent to $\lambda_{\Gamma}$ (because $\Gamma$
is C$^*$-simple).  
This is an illustration of Glimm's theorem 
(see Part VI of the appendix). 

   For $t \in ]0,1[$ and the corresponding irreducible representation
$\pi^{Co}_t$ of  the complementary series of $G$, the restriction 
$\pi^{Co}_t \vert \Gamma$ is not weakly equivalent to $\lambda_{\Gamma}$. Moreover,
if $t < t'$,  the representation $\pi^{Co}_t \vert \Gamma$ is not weakly contained in
$\pi^{Co}_{t'} \vert \Gamma$. \smallskip

   For all this, see \cite{BeHa--94} and \cite{Bek--97}.

\medskip

\proclaim{23.~Example} 
Other irreducible representations 
weakly equivalent to a left regular repre\-sentation.
\endproclaim

   Consider more generally a non-compact semi-simple connected real Lie group $G$ 
with Iwasawa decomposition $G = KAN$, 
a maximal amenable subgroup $B = MAN$ of $G$, 
the characters $\chi_s$ of $A$ viewed as one-dimensional representations of $B$, 
and the induced unitary representations $\rho = \operatorname{Ind}_B^G(\chi_s)$; 
restrict attention to those
$\chi_s$ which provide irreducible $\rho$. 
Then the restrictions $\rho \vert \Gamma$ are 
irreducible, pairwise non-equivalent \cite{CoSt-91}, 
and all are weakly equivalent to $\lambda_{\Gamma}$. 

\medskip

The next example is a remarkable result of Yoshizawa:
see \cite{Yos--51}, as well as Theorem VII.6.5 in \cite{Dav--96}.

\proclaim{24.~Example} 
For a non-abelian free group, 
the universal representation $\pi_{uni}$ of Section~2
is weakly equivalent to an irreducible representation.
\endproclaim

   With the method of \cite{Dav--96}, 
it is not difficult to extend this to a few other cases, 
for example to groups of the form
$\Gamma_0 \ast \bold Z \ast (\bold Z/2\bold Z)$ or
$\Gamma_0 \ast \bold Z \ast \bold Z$
(where $\Gamma_0$ is an arbitrary countable group,
and where $\ast$ indicates a free product).
Observe also that, if $\Gamma$ is an amenable icc group,
its universal representation 
is weakly equivalent to its regular representation,
and therefore weakly equivalent to an irreducible represenation
(Proposition 21).
However, we know very little in general:

\proclaim{25.~Problem} 
Find other groups for which the universal representation
is   weakly equivalent to an irreducible one. 
\endproclaim

\bigskip

\subhead
III.~On diagonal and arbitrary matrix coefficients
\endsubhead

\medskip

It is important that the definition of weak containment involves    
{\it diagonal matrix} coefficients, and not arbitrary matrix coefficients. 
To show the difference, let us consider the left-regular
representation $\lambda_{\Gamma}$ of a non-amenable finitely-generated group 
$\Gamma$. 

   Let $A(\Gamma)$ be the set of {\it all coefficients} of the form  
$$
\varphi_{\xi,\eta} \, : \, \left\{
\aligned \Gamma \ &\longrightarrow \hskip1truecm \bold C \\
       \gamma \ &\longmapsto \ 
    \langle \xi \mid \lambda_{\Gamma} (\gamma) \eta \rangle   
\endaligned \right.
$$
with $\xi,\eta \in \ell ^2(\Gamma)$;
if $\tilde \eta$ is defined by 
$\tilde \eta (\gamma) = \overline{\eta(\gamma^{-1})}$,
observe that $\varphi_{\xi,\eta}$
is the convolution $\xi \ast \tilde \eta$.
It is a standard result that this {\it set} is a {\it linear subspace} 
of the commutative C$^*$-algebra 
$\Cal C_0(\Gamma)$ of all complex-valued functions on $\Gamma$ 
vanishing at infinity (Page 218 in \cite{Eym--64}).

On the one hand, since $\lambda_{\Gamma} \otimes \lambda_{\Gamma}$ 
is unitarily equivalent to a multiple of $\lambda_{\Gamma}$ 
(see e.g.\ Complement~13.11.3 in \cite{DC$^*$--69}),  
it is a consequence of the Stone-Weierstrass theorem that
$A(\Gamma)$ is a dense  subalgebra of $\Cal C_0(\Gamma)$.
[It is the so-called {\it Fourier algebra} of Eymard.]

   On the other hand, since $\Gamma$ is not amenable, 
it follows from a classical  result  of Kesten \cite{Kes--59} 
that there is no sequence of {\it diagonal} matrix coefficients, 
namely no sequence of coefficients of the form 
$\left( \varphi_{\xi_n,\xi_n} \right)_{n \ge 1}$, which can converge towards a
constant function, uniformly on all finite subsets of $\Gamma$. 
Indeed, for $S$ a finite set of generators of $\Gamma$ such that $S^{-1} = S$,
say of size $\vert S \vert = d$, Kesten has shown that the self-adjoint
operator $h = \frac{1}{d} \sum_{s \in S} \lambda_{\Gamma}(s)$ has
spectrum contained in $[-1,1-\epsilon]$ for some $\epsilon > 0$, or in other
words that
$$
   \frac{1}{d} \sum_{s \in S} \varphi_{\xi,\xi} (s) \, = \, \frac{1}{d}
   \sum_{s \in S} \left\langle \xi \vert \lambda_{\Gamma}(s) \xi \right\rangle
   \, \le \, 1 - \epsilon
$$
for any unit vector $\xi \in \Cal H$. Hence there exists $s \in S$
with $\Re \left( \varphi_{\xi,\xi} (s) \right) \le 1 - \epsilon$. 
As $\varphi_{\xi,\xi} (1) = 1$, this shows that the restriction of 
$\varphi_{\xi,\xi}$ to the finite set $S \cup \{1\}$ is bounded away from
a constant function.

\bigskip

\subhead
IV.~On a definition of Zimmer
\endsubhead

\medskip

Definition 7.3.5 of \cite{Zim--84} offers a {\it different definition} of
weak containment. Consider two unitary representations $\rho,\sigma$ of some
group $G$. First, to a finite orthonormal set $\xi_1 , \hdots , \xi_n$ of 
$\Cal H_{\rho}$, we associate a continuous function 
$$
\varphi_{(\xi_1,\hdots,\xi_n)} \, : \, \left\{
\aligned G \ &\longrightarrow \hskip1truecm Mat_n(\bold C) \\
         g \ &\longmapsto \ \left( \langle \xi_j \mid \rho(g) \xi_k \rangle
                            \right)_{1 \le j,k \le n}
\endaligned \right.
$$
which is called a {\it $n$-by-$n$ submatrix of $\rho$.}
Define then $\rho$ to be 
{\it weakly contained in $\sigma$ in the sense of Zimmer,} and write
$\rho \overset{Z} \to \prec \sigma$, 
if any submatrix of $\rho$ is a limit of submatrices of
$\sigma$, uniformly on compact subsets of $G$. Two unitary representations are
{\it weakly equivalent in the sense of Zimmer,} 
written $\rho \overset{Z} \to \sim \sigma$, if 
$\rho \overset{Z} \to \prec \sigma$ and
$\sigma \overset{Z} \to \prec \rho$. \par

   The definition of our Section 1 is that of Fell 
(see \cite{Fel--60} and \cite{Fel--64}) and
Dixmier (see Items 3.4.5 and 18.1.4 in \cite{DC$^*$--69}). \par

   The two notions are distinct. For example, if $\sigma$ and $\rho$ are
finite-dimensional, then $\rho \overset{Z} \to \prec \sigma$ \`a la Zimmer if and
only if $\rho$ is a subrepresentation of $\sigma$ in the naive sense, whereas
$\rho \oplus \rho \prec \rho$.  It is a fact that the
two notions are equivalent for irreducible representations, and indeed 
$\rho \prec \sigma \Longleftrightarrow \rho \overset{Z} \to \prec \sigma$ whenever 
$\rho$ is irreducible. It is another fact that $\rho \prec \sigma$ \`a la
Fell-Dixmier if and only if  $\rho \overset{Z} \to \prec \infty \, \sigma$ \`a la
Zimmer  (where $\infty \, \sigma$ stands for an
orthogonal sum of infinitely many copies of the representation
$\sigma$).

\bigskip

\subhead 
V.~Other notions of equivalence
\endsubhead 

\medskip

      In case they are both {\it irreducible,} two unitary representations
$\rho,\sigma$ of a second-countable locally compact group $G$ which are weakly
equivalent are also {\it approximately unitarily equivalent relative to the
compact operators,} in the following sense~: there exists a sequence of unitary
operators $$
   \left( U_n : \Cal H_{\rho} \longrightarrow \Cal H_{\sigma} \right)_{n \ge 1}
$$
such that
$$
\aligned
   &U_n \rho(x) U_n^* \, - \, \sigma(x) \quad \text{is compact} \quad
     \forall \ n \ge 1 \quad \text{and} \quad \forall \ x \in C^*_{max}(G) ,\\
   &\lim_{n \to \infty} \norm{U_n \rho(x) U_n^* \, - \, \sigma(x) }
      \, = \, 0 \quad \forall \ x \in C^*_{max}(G) .
\endaligned
$$
This follows from a deep result of D. Voiculescu \cite{Voi--76}; see also
Proposition VII.3.4 and Theorem II.5.6 in \cite{Dav--96}.

\medskip 

   Note that there exist at least four notions of equivalence for two unitary
representations $\rho,\sigma$ of a second countable locally compact group $G$~:

\smallskip

   (i) \phantom{iiv} unitary equivalence $\rho \approx_u \sigma$, \par
   (ii) \phantom{iv} weak equivalence $\rho \sim \sigma$, \par
   (iii) \phantom{v} approximate unitary equivalence relative to the compact
       $\rho \sim_{\Cal K} \sigma$, \par
   (iv) \phantom{ii} quasi-equivalence $\rho \approx_{q.e.} \sigma$
       (see \S \ V in \cite{DC$^*$--69}). 

\smallskip\noindent In the particular case $\rho$ and $\sigma$ are irreducible,
we have 
$$
\aligned
   \rho \, \approx_u \, \sigma  &\qquad \Longleftrightarrow \qquad
     \rho \, \approx_{q.e.} \, \sigma \hskip1truecm
     \text{(Proposition 5.33 in \cite{DC$^*$--69}),} \\ 
   \rho \, \sim \, \sigma  &\qquad  \Longleftrightarrow \qquad 
     \rho \, \sim_{\Cal K} \, \sigma \hskip1.25truecm
     \text{(Voiculescu's result recalled above),} \\ 
   \rho \, \approx_u \, \sigma &\qquad \Longrightarrow \qquad
     \rho \, \sim \, \sigma \hskip1.6truecm
     \text{(obvious),} \\
   \rho \, \approx_u \, \sigma \ 
            &\overset{ \text{for  } G \text{  of type I} } \to \Longleftarrow  \ \,
     \rho \, \sim \, \sigma \hskip1.6truecm 
     \text{(see Part VI of the appendix).}
\endaligned
$$

\bigskip

\subhead 
VI.~The dual of a group
\endsubhead 

\medskip

   Let $A$ be a C$^*$-algebra, $\hat A$ the set of equivalence classes of
irreducible representations of $A$ and $Prim(A)$ the set of primitive two-sided
ideals in $A$. As primitve ideals are precisely the kernels of irreducible
representations, there is a canonical mapping from $\hat A$ onto $Prim(A)$. 
On $Prim(A)$, there is a natural topology called the {\it hull-kernel} or the
{\it Fell-Jacobson} topo\-logy, for which the closure of a subset $S$ is defined
to be the set of primitive ideals of $A$ containing $\cap_{\Cal I \in S} \Cal I$.
By definition, the {\it spectrum} of $A$ is the set $\hat A$ together with the
topology pulled back from the hull-kernel topology on the primitive ideal space.
There are other definitions of the same topology; moreover, if $A$ is separable,
this topology is compatible with the appropriate  \lq\lq Mackey Borel
structure\rq\rq \ on $\hat A$. (See Theorems 3.4.11 and 3.5.8, as well as
Proposition 3.1.3, in \cite{DC$^*$--69}.) This topology is useful to decompose
arbitrary representations as \lq\lq direct integrals\rq\rq \ of irreducible ones,
and also for the study of the structure of C$^*$-algebras. \par

   The dual $\hat G$ of a second-countable locally compact group $G$ is the
spectrum of its maxi\-mal C$^*$-algebra $C^*_{max}(G)$. For irreducible
representations $\rho,\sigma$ of $G$, it follows from the definitions
that $\rho \prec \sigma$ if and only if $\rho \in \overline{ \{\sigma\} }$.
\par

   The following theorem, for its main part due to Glimm, is fundamental.
See \cite{Gli--61}, \S \ 9 of \cite{DC$^*$--69}, Theorem 6.8.7 in \cite{Ped--79},
and \cite{Mar--75} (of which the last statement of the theorem below is a special
case). In the particular case of the maximal C$^*$-algebra of a second-countable
locally compact group, we can replace \lq\lq irreducible representation of
$A$\rq\rq \ by \lq\lq irreducible unitary representation of $G$\rq\rq \
and $\hat A$ by $\hat G$.

\proclaim{26.~Theorem} 
Let $A$ be a separable C$^*$-algebra. The following properties
are equivalent~: 
\roster
\item"(i)"   
   the canonical mapping  $\hat A \longrightarrow Prim(A)$ is a  bijection, 
\item"(ii)"
   the image of any irreducible representation of $A$ contains some 
   non-zero compact operator, 
\item"(iii)"
   the image of any irreducible representation of $A$ contains all 
   compact operators, 
\item"(iv)"
   the appropriate Borel structure on $\hat A$ is  countably  separated, 
\item"(v)"
   two irreducible representations of $A$ are equivalent if and only if 
   they are weakly equivalent, 
\endroster
and they imply 
\roster
\item"(vi)"
   the spectrum $\hat A$ contains an open dense subset which is locally compact.
\endroster
Moreover, if $A$ does {\rm not} have these properties, for any injective factor
$M$ of type $II_{\infty}$ or $III$, there exists an irreducible representation
$\rho$ of $A$ such that $\rho(A)$ generates $M$. 
\endproclaim

\bigskip

\subhead 
VII.~Groups of Type I
\endsubhead 

\medskip
   
   A C$^*$-algebra which has the properties of Theorem~26 is said to be
{\it of type I.} A second countable locally compact group $G$ is of type I if
$C^*_{max}(G)$ is of type I. \par

   Groups of type I include include 
connected real Lie groups which are either nilpotent  or real algebraic, 
the \lq\lq$ax+b$-group\rq\rq \ (the  affine group over $\bold R$) 
and various other soluble real Lie groups (caracterized in \cite{AuKo--71}),  
$p$-adic linear groups which are either reductive \cite{Ber--74} 
or soluble \cite{GoKa--79}, 
large classes of adelic groups \cite{GGP--69}, 
countable groups which have normal abelian groups of finite index, 
but {\it no other countable group} \cite{Tho--68}. 

There is an existence theorem according to which 
we can \lq\lq decompose\rq\rq \ any unitary repre\-sentation in a separable
Hilbert space as a direct integral of irreducible representations. 
For a group of type I, this decomposition is essentially unique 
(Sections 8.5 and 8.6 in \cite{DC$^*$--69}).
But this uniqueness fails completely for other groups;
specific examples are worked out
for the affine group over $\bold Q$ in Section 3.5 of  \cite{Mac--76}
and for the free group on two generators in Chapter 19 of \cite{Rob--83}.\par

    A C$^*$-simple group $\ne \{1\}$ is never of type I.

\bigskip

\subhead 
VIII.~The reduced dual 
\endsubhead 

\medskip

    The {\it reduced dual} of a second-countable locally compact group $G$ is the
spectrum $\hat G_{red}$ of its reduced C$^*$-algebra,
or equivalently the space of those irreducible unitary representations of $G$
which are weakly contained in $\lambda_G$. 
The canonical embedding of $\hat G_{red}$ in $\hat G$  has a closed image, 
and is onto if and only if $G$ is amenable. \par
 
   When $\Gamma$ is C$^*$-simple, any irreducible representation in 
$\hat \Gamma_{red}$ is weakly equivalent to $\lambda_{\Gamma}$, 
and $\hat \Gamma_{red}$ is the closure of any of  its points ! 
See Examples 22 and 23 for specific  examples. \par

   For irreducible representations of free groups 
which are weakly equivalent to the  regular representation, see 
\cite{PySz--86}, \cite{Szw--88}, \cite{KuSt--92}, and
\cite{BuHa--97, Proposition 4.1}. 

   There is  a criterium, involving an appropriate growth condition, 
for a function of positive type on a non-abelian free group 
to be associated with the reduced dual of the group \cite{Haa--79}. 
The result and its proof hold for some other groups, 
such as Gromov-hyperbolic Coxeter groups, 
as observed in Section 2 of \cite{JoVa--91}.

\bigskip

\subhead 
IX.~The Dixmier property 
\endsubhead 

\medskip

   Let $A$ be a C$^*$-algebra with unit $1_A$. For $a \in A$, set
$$
   C_A(a) \, = \, \overline{\text{conv}} 
   \big\{ uau^* \mid u \in A , \, uu^* = u^*u = 1_A \big\} \, \cap \,
   \bold C 1_A
$$
where $\overline{\text{conv}}$ indicates the norm-closed convex hull.
Then $A$ has the {\it Dixmier property} if 
$C_A(a) \ne \emptyset \ \forall \ a \in A $.
Dixmier has shown that a C$^*$-algebra which has a tracial state and which has
the Dixmier property has a unique tracial state and is simple 
(see \cite{DvN--69}, Chapter III, Sections 5.1, 5.2, 8.5 and 8.6).
A result of \cite{HaZs--84} shows that, if $A$ has a unique tracial state $tr$
(e.g.\ if $A$ is the reduced C$^*$-algebras of a Powers group), 
then $C_A(a)$ is reduced to a unique element, which is $tr(a)1_A$. \par

   The main step in the proof of Theorem 14 shows that the set $C_A(U)$,
which is obviously inside the ideal $\Cal I$, contains an invertible element
($U$ is the element of the ideal $\Cal I$ chosen at the beginning of the proof
of Theorem 14).
As 
$$
   tr(U) \, = \,  tr\bigg(1 + \sum_{x \in \Gamma, x \ne 1} z_x \lambda(x) \bigg) 
         \, = \, 1 ,
$$
we have indeed $C_A(U) = \{1_A\}$, and the strategy of the proof of Theorem 14
should not be surprising. \par

   S. Popa has used his \lq\lq relative Dixmier property\rq\rq \ 
to show Claims $(iii)$ and $(iv)$ of our Proposition 19.       

\bigskip

\subhead 
X.~On the icc condition
\endsubhead 

\medskip

   For a group $\Gamma$, denote by $\Gamma_f$ the union of the finite conjugacy
classes of $\Gamma$. It is clearly a subgroup of $\Gamma$ which is normal,
indeed characteristic; let us check that $\Gamma_f$ is also amenable.

   It is enough to check that any subgroup $\Delta$ of $\Gamma_f$
generated by a finite set $\{s_1,\hdots,s_n\}$ is amenable.
Since the conjugacy class of each $s_j$ is finite
(both in $\Gamma$ and in $\Delta$),
the centralizer $Z_{\Delta}(s_j)$ of $s_j$ in $\Delta$
is a subgroup of finite index;
thus the intersection $Z = \bigcap_{1 \le j \le n}Z_{\Delta}(s_j)$
is also of finite index in $\Delta$.
Since $Z$ is contained in the centre of $\Delta$,
the group $\Delta$ contains a normal abelian subgroup of finite index,
and is therefore amenable.

  Recall from the footnote to Item~13 that $\Gamma$ is icc
if $\Gamma \ne \{1\}$ and $\Gamma_f = \{1\}$.
It follows from Proposition~3 that a C$^*$-simple group is icc.
The converse is far from being true; for example, the group $\Sigma$
of permutations of finite support of an infinite countable set
is amenable (because it is locally finite), and therefore not C$^*$-simple,
though $\Sigma$ is clearly an icc group.

   The Murray and von Neumann lemma from \cite{ROIV} already quoted
in Part~II of the appendix shows that the von Neumann algebra $W^*(\Gamma)$
is a factor 
\footnote{
Recall that a von Neumann algebra $M$ is a factor of type $II_1$
if it is infinite dimensional, simple, and finite in the following sense:
for $x,y \in M$, we have $xy = 1$ if and only if $yx = 1$.
}
of type $II_1$ if and only if the group $\Gamma$ is icc.
The previous observations show that the icc condition is also necessary
for the reduced C$^*$-algebra $C^*_{red}(\Gamma)$ to be simple.
Short of suggesting general necessary and sufficient conditions on $\Gamma$
for $C^*_{red}(\Gamma)$ to be simple, we would like to end this report
by stating three problems.

\medskip

\proclaim{27.~Problem}
Let $\Gamma$ be the free product with amalgamation defined by 
two groups $\Gamma_1,\Gamma_2$ and a common subgroup $\Gamma_0$.
Find necessary and
\footnote{
Short of this, find necessary {\it or} sufficient conditions!
}
sufficient conditions for $\Gamma$ to be icc,
and for $\Gamma$ to be C$^*$-simple.
\endproclaim

For the particular case of free products ($\Gamma_0 = \{1\}$), 
here is one formulation of the classical solution of Problem 27.
Assume that neither $\Gamma_1$ nor $\Gamma_2$ is reduced to one element;
for the free product $\Gamma = \Gamma_1 \ast \Gamma_2$, 
the three following conditions are equivalent:
\roster
\item"(i)"
  $\Gamma$ is not the infinite dihedral group 
  $(\bold Z /2\bold Z) \ast (\bold Z /2\bold Z)$;
\item"(ii)"
   $\Gamma$ is icc;
\item"(iii)"
   $\Gamma$ is C$^*$-simple.
\endroster
(This follows from \cite{PaSa--79}, or from Corollary 12.i.)

   For the case of Problem~27, an almost obvious sufficient condition
for $\Gamma$ to be icc
is that one at least of the two groups $\Gamma_1,\Gamma_2$ is icc.
There is a sufficient condition for $\Gamma$ to be C$^*$-simple
in \cite{B\'ed--84}.

\medskip

\proclaim{28.~Problem}
Let $\Gamma$ be the HNN extension defined by 
a group $\Gamma_1$
and an isomorphism $\varphi$ from a subgroup $\Gamma_0$ 
onto another subgroup $\Gamma'_0$  of $\Gamma_1$.
Find necessary and
sufficient conditions for $\Gamma$ to be icc,
and for $\Gamma$ to be C$^*$-simple.
\endproclaim

There are partial answers for the icc condition in \cite{Sta} and \cite{HaPr}.

\medskip

Let $M$ be an orientable connected compact manifold of dimension $3$
and let $\Gamma$ denote its fundamental group; 
we assume that $\Gamma$ is infinite.
In case $M$ is a Seifert manifold, it is known that 
$\Gamma$ contains an infinite cyclic normal subgroup,
so that in particular $\Gamma$ is not icc.
We have shown  that, conversely,
if $\Gamma$ is not the fundamental group of a Seifert $3$-manifold,
then $\Gamma$ is icc
(see \cite{HaPr}, which contains also a discussion of the non-orientable case).

\proclaim{29.~Problem}
Let $\Gamma$ be the fundamenal group of 
an orientable connected compact manifold of dimension $3$ 
which is not the fundamental group of a Seifert $3$-manifold.
When is $\Gamma$ C$^*$-simple?
\endproclaim

\bigskip

   I am grateful to 
Bachir~Bekka for numerous comments on the material of this exposition,
and to Alain~Valette for useful observations.

\enddocument

\vskip1cm

\subhead
Pas encore fait le 1er septembre
\endsubhead

\medskip
\noindent $\circ$
Caract\'erisation des groupes $\Gamma$ dont le dual $\hat \Gamma$
est d\'enombrable ??? (pour VIII).

\medskip
\noindent $\circ$
Revoir les r\'ef\'erences de l'appendice VIII.
(Gard\'e toutes le 7/9.)

\enddocument